\documentclass[twocolumn]{autart} 

\usepackage{savesym}
\savesymbol{AND}
\usepackage{graphicx}
\usepackage{amsmath}
\usepackage{amssymb}
\usepackage{subfig}
\usepackage{epstopdf}
\usepackage{algorithmic}
\usepackage{algorithm}
\usepackage{bm}
\usepackage{cite}

\newtheorem{theorem}{Theorem}
\newtheorem{lemma}{Lemma}
\newtheorem{assumption}{Assumption}
\newtheorem{definition}{Definition}
\newtheorem{corollary}{Corollary}
\newtheorem{proposition}{Proposition}
\newtheorem{remark}{Remark}

\DeclareMathOperator*{\argmin}{arg\,min}
\DeclareMathOperator{\proj}{proj}

\newcommand{\RR}{\mathbb{R}}
\newcommand{\xx}{\mathbf{x}}
\newcommand{\oxx}{\overline{\xx}}
\newcommand{\gt}{\mathbf{g}_{t}}
\newcommand{\ee}{\mathbf{e}}
\newcommand{\XX}{\mathcal{X}}
\newcommand{\reg}{\texttt{Reg}}
\newcommand{\Tr}{^{\top}}
\newcommand{\xxt}{\xx_{t+1}}
\newcommand{\oxt}{\oxx_{t+1}}
\newcommand{\dt}{\mathbf{d}_{t+1}}
\renewcommand{\ss}{\mathbf{s}}
\newcommand{\cc}{\mathbf{c}}

\newcommand*{\QED}{\hfill\ensuremath{\square}}

\begin{document}

\begin{frontmatter}
%\runtitle{Insert a suggested running title}  % Running title for regular 
                                              % papers but only if the title  
                                              % is over 5 words. Running title 
                                              % is not shown in output.

\title{Predictive Online Convex Optimization\thanksref{footnoteinfo}}

\thanks[footnoteinfo]{This work was partly done while A. Lesage-Landry was visiting The University of Melbourne, Australia. This work was funded by the Fonds de recherche du Qu\'ebec -- Nature et technologies, the Ontario Ministry of Research, Innovation and Science and the Natural Sciences and Engineering Research Council of Canada. This paper was not presented at any IFAC 
meeting. Corresponding author A. Lesage-Landry Tel.: +1 416-978-6842; fax: +1 416-978-1145 }

\author[ucb]{Antoine Lesage-Landry}\ead{alesagelandry@berkeley.edu},    % Add the 
\author[is]{Iman Shames}\ead{iman.shames@unimelb.edu.au},               % e-mail address 
\author[ut]{Joshua A. Taylor}\ead{josh.taylor@utoronto.ca}  % (ead) as shown

\address[ucb]{Energy \& Resources Group, University of California, Berkeley, USA}                                     
\address[is]{Department of Electrical and Electronic Engineering, University of Melbourne, Parkville, Australia}             % full addresses
\address[ut]{The Edward S. Rogers Sr. Department of Electrical and Computer Engineering, University of Toronto, Toronto, Canada}  % Please supply     


\begin{abstract}                          
We incorporate future information in the form of the estimated value of future gradients in online convex optimization. This is motivated by demand response in power systems, where forecasts about the current round, e.g., the weather or the loads' behavior, can be used to improve on predictions made with only past observations. Specifically, we introduce an additional predictive step that follows the standard online convex optimization step when certain conditions on the estimated gradient and descent direction are met. We show that under these conditions and without any assumptions on the predictability of the environment, the predictive update strictly improves on the performance of the standard update. We give two types of predictive update for various family of loss functions. We provide a regret bound for each of our predictive online convex optimization algorithms. Finally, we apply our framework to an example based on demand response which demonstrates its superior performance to a standard online convex optimization algorithm.
\end{abstract}

\end{frontmatter}

\section{Introduction}
Online convex optimization (OCO) has found applications in fields like network resource allocation~\cite{chen2017online,chen2018bandit,cao2018virtual} and demand response in power systems~\cite{kim2017online,lesage2018setpoint}. It is used for sequential decision-making when contextual information or feedback is only revealed to the decision maker at the end of the current round. Theoretical results showing that OCO algorithms have bounded regret guarantee the performance of these algorithms under mild assumptions. 

In many applications, the decision maker has access to both revealed past information and estimates about future rounds. For example, in power systems, weather forecasts or historical load patterns can be used to estimate the future regulation needs~\cite{mathieu2013state,callaway2009tapping}. In this work, we present the predictive online convex optimization (POCO) framework. POCO works under the assumption that an estimate of the gradient of the loss function for the next round is available to the decision maker. In POCO, a standard OCO update is first applied using past information to compute the next decision. Then, the decision maker checks the quality of the estimated information available to them. If the estimated gradient is considered accurate enough, the decision maker implements an additional projected gradient step based on the estimated gradient to improve their decision for this round. This last step is referred as the predictive update.

We introduce explicit criteria for determining if the quality of the estimated gradient is high enough to guarantee an improvement over a standard OCO step when the predictive update is applied. A regret bound is obtained for all our algorithms. We conclude this work by presenting numerical examples where a POCO algorithm is used to improve on the performance of demand response with standard OCO. This example is motivated by the fact that a load aggregator often has access to an estimate of the power imbalance they have to counteract for regulation purposes.

\textit{Literature review.} Recent work in online convex optimization has focused on including prior or future information. Reference~\cite{rakhlin2013optimization}, which builds on~\cite{chiang2012online}, assumes that the problem's unknown and uncertain parameters follow a predictable process plus some noise~\cite{rakhlin2013online} for their OCO algorithm. As in our setting, a second update with an estimated gradient-like term follows a mirror descent update. This second update is used by the algorithm in every step regardless of the quality of the estimated gradient. For this reason, the algorithm is referred to as optimistic. Optimistic algorithms were also studied in~\cite{steinhardt2014adaptivity,mohri2016accelerating,yang2016optimistic}. No conditions are provided about the estimated gradient in this case except that it comes from past observations and/or side information via an oracle. The authors of~\cite{rakhlin2013optimization} show that the optimistic mirror descent can lead to a tighter bound than a standard online mirror descent algorithm if the process is indeed predictable. In~\cite{jadbabaie2015online}, the authors provide a dynamic regret bound for the optimistic mirror descent. There is, however, no guarantee that in a given round the optimistic update does not do worse than the standard OCO update. An algorithm similar to~\cite{rakhlin2013optimization} is given in~\cite{ho2016accelerating}. In their work, they make the stronger assumption in which the exact gradient of the next round loss function is available and then provide a static regret bound for their setting. This differs from our setting in that we provide dynamic regret-bounded algorithms and use an estimated gradient which entails less restrictive assumptions. Several other authors have studied different ways to incorporate future information in OCO like using state information~\cite{hazan2007online} or the direction of the loss function's gradient in an online linear optimization setting~\cite{dekel2017online}.

The projected gradient descent, inexact gradient descent, and proximal algorithms~\cite{bertsekas1997nonlinear,schmidt2011convergence,bertsekas2015convex} from conventional convex optimization resemble our setting. These algorithms differ from ours because they aim to minimize the same objective function throughout all descent steps. In OCO, we minimize a sequence of objective functions $\left\{f_t\right\}_{t=1}^T$ and at each time $t$ provide a decision to minimize the current loss function. The loss function in a given round is only observed after we have committed to a decision. OCO will be introduced formally in Section~\ref{sec:background}.

Model predictive control (MPC)~\cite{garcia1989model,borrelli2017predictive} is another widely-used sequential decision-making framework. In MPC, the decision maker solves to optimality a receding horizon optimization problem that relies on models of future round loss functions. This thus requires significantly more contextual information and computational resources. These limitations are absent in OCO, making it a more suitable tool for real-time decision making with small computational resources. 

Because we characterize conditions under which the predictive step improves performance, we guarantee improvement over conventional OCO and require no predictability assumptions. These conditions can be checked at each round of OCO, and if satisfied, the predictive update is implemented. In sum, in this work we make the following contributions:
\begin{itemize}
  \item We introduce a novel predictive online convex optimization framework and provide conditions for when to use side information.
  \item We propose a predictive update with a predetermined step size for loss functions that have a Lipschitz gradient. We show that this update leads to a strict improvement over an OCO update when used (Section~\ref{sec:predetermined}).
  \item We give a predictive update with backtracking line search that applies to a broader family of problems. We show that it leads to strict improvement over an OCO update (Section~\ref{sec:backtracking}).
  \item We obtain sublinear regret bounds in the number of rounds for all algorithms.
  \item We apply our framework to demand response in power systems and find that it outperforms a standard OCO algorithm (Section~\ref{sec:example}).
\end{itemize}

\section{Background}
\label{sec:background}
In OCO, one must make a decision at each round to minimize their cumulative loss~\cite{shalev2012online,hazan2016introduction}. The current round's loss function and any other round-dependent parameters are not available at the moment when the decision is made. Only information about previous rounds can be used to make the decision. Once the decision has been made, information about the current round is observed.

Let $t$ denote the current round index and $T$ be the time horizon. Let $\XX \subset \mathbb{R}^N$, $N \in \mathbb{N}$, be the decision set, and let $\xx_t \in \XX$ be the decision variable at time $t$. We denote the differentiable convex loss function by $f_t(\xx_t)$ for $t=1,2,\ldots, T$. Let $\| \cdot \|$ be the Euclidean norm. We denote the projection operator onto the set $\mathcal{Y}$ as $\proj_\mathcal{Y} ( \xx ) \in \argmin_{\mathbf{y} \in \mathcal{Y}} \left\| \xx - \mathbf{y} \right\|$. 

The goal of the decision maker is to sequentially solve the following sequence of problems:
\begin{equation}
\min_{\xx_t \in \XX} f_t(\xx_t)
\label{eq:oco_problem}
\end{equation}
for $t=1,2,\ldots, T$. The decision maker observes the loss function $f_t$ after choosing $\xx_t$. For this reason, even if the loss function has a simple form, an analytical solution to the round optimization problem~\eqref{eq:oco_problem} is not obtainable. The decision $\xx_t$ is computed using a gradient descent-based~\cite{zinkevich2003online}, mirrored descent-based~\cite{duchi2010composite} or Newton step-based rule~\cite{hazan2007logarithmic}. For example, in the online gradient descent (\texttt{OGD})~\cite{zinkevich2003online} algorithm, the decision at round $t+1$, $\xx_{t+1}$, is given by the update:
\begin{equation}
\xx_{t+1} = \proj_{\XX} \left(\xx_t - \eta \nabla f_t(\xx_t)\right), \label{eq:ogd_update}
\end{equation}
where $\eta \propto T^{-1/2}$ to guarantee a sublinear upper bound on the dynamic regret of \texttt{OGD}~\cite[Theorem 2]{zinkevich2003online}.

Throughout this work, we make the following assumptions~\cite{zinkevich2003online,shalev2012online,hazan2016introduction}.
\vspace{-0.2cm}
\begin{assumption}
The set $\XX$ is convex and compact. \label{ass:compact}
\label{assumption:compact}
\end{assumption}
\vspace{-0.2cm}
The decision set $\XX$ represents all constraints on $\xx_t$. In this version of OCO, we only consider time-invariant constraints.
\vspace{-0.2cm}
\begin{assumption}
The loss function is $B$-bounded: $|f_t(\xx_t)| \leq B$ for $t=1,2,\ldots, T$ and $B<\infty$.
\label{assumption:bounded}
\end{assumption}
\vspace{-0.2cm}

\begin{assumption}
The gradient of the loss function is $G$-bounded: $\|\nabla f_t(\xx_t)\| \leq G$ for $t=1,2,\ldots, T$ and $G<\infty$.
\label{assumption:grad_B}
\end{assumption}
\vspace{-0.2cm}

As a consequence of Assumption~\ref{ass:compact}, the decision variable is also $X$-bounded: $\| \xx_t \| \leq X$ for $t=1,2,\ldots, T$. We define the diameter of the compact set $\XX$ as $\mathrm{diam} \  \XX  = \sup\Big\{ \|\xx - \mathbf{y} \| \Big| \xx, \mathbf{y} \in \XX \Big\}$and let $D = \mathrm{diam} \  \XX$, a positive scalar. The remainder of the assumptions will be stated when a specific technical result requires it.

The design tool of OCO algorithms is the regret~\cite{shalev2012online,hazan2016introduction}. In this work, we use the dynamic regret~\cite{chen2017online,jadbabaie2015online,zinkevich2003online,mokhtari2016online}:
\begin{equation}
\texttt{Reg}^d_T = \sum_{t=1}^T f_t(\xx_t) - f_t(\xx_t^\ast), \label{eq:regret_d}
\end{equation}
where $\xx_t^\ast \in \argmin_{\xx \in \XX} f_t(\xx)$. The dynamic regret compares the loss suffered by the decision maker to optimal performance in each round. Other versions of the regret exists, e.g., static regret~\cite{shalev2012online,hazan2016introduction,zinkevich2003online}, which is defined in terms of the optimal stationary decision, $\xx^\ast \in \argmin_{\xx \in \XX} \sum_{t=1}^T f_t(\xx)$ in~\eqref{eq:regret_d}. In this work, we only consider the dynamic regret because it yields a stronger theoretical guarantee. This theoretical guarantee is also more relevant in the context of time-varying optimization. For this reason, we refer to the dynamic regret, $\reg^d_T$, simply as the regret. Note that a bounded dynamic regret implies a bounded static regret~\cite{chen2017online}. The goal when designing an OCO algorithm is to show that the regret is sublinearly bounded above in the number of rounds. An OCO algorithm with a sublinearly bounded dynamic regret in the number of rounds will on average perform as well as the round optimal decision at each round~\cite{shalev2012online,hazan2016introduction,hazan2007logarithmic}. 

We conclude this section by defining the quantity $V_T = \sum_{t=2}^T \left\| \xx_t^\ast - \xx_{t-1}^\ast \right\|$. The term $V_T$ quantifies the variation of the optimal predictions through all rounds.

\section{Predictive OCO}
We now introduce our POCO framework. We let $\oxt \in \XX$ be the decision computed by an OCO algorithm in round $t$. This OCO algorithm can be, for example, the aforementioned \texttt{OGD}. The decision $\oxt$ is then given by the update~\eqref{eq:ogd_update}. In POCO, we consider an $\epsilon$-forecaster introduced in Assumption~\ref{def:eps_fore}. Let $\gt(\oxt) \in \RR^N$ be the estimated gradient of the loss function $f_{t+1}$ at $\oxt$.

\begin{assumption}[$\epsilon$-forecaster]
The $\epsilon$-forecaster has access to an estimate of the gradient of the next round's loss function evaluated at $\oxt$, $\gt\left(\oxt\right)$, and the maximum estimation error, $\epsilon > 0$, such that $\|\gt\left(\oxt\right) - \nabla f_{t+1}(\oxt)\|\leq \epsilon$ for all time $t = 1, 2, \ldots, T$.
\label{def:eps_fore}
\end{assumption}
In other words, we consider a forecaster that has access to limited information about the next round in the form of $\gt\left(\oxt\right)$. This could represent a prediction based on historical data, e.g., a weather or demand forecast. The decision maker uses this information to improve on the OCO update. For conciseness, we denote the estimated gradient by $\gt$. We omit its dependency on $\oxt$ because it is always evaluated at the OCO update output, $\oxt$, and no other points. The decision maker can meet Assumption~\ref{def:eps_fore} by relying on an exogenous model to estimate the gradient $\nabla f_{t+1}\left( \oxt \right)$. In the context of demand response, historical data of the load's consumption and generator output’s patterns, weather history and the historical values of the gradient, for example, can be used to build a statistical model to estimate the value of the $\nabla f_{t+1}$ at the decision given by OCO update. The parameter $\epsilon$ can then be set according to, for example, a high confidence interval or a worst-case performance parameter. The forecaster would then provide $\gt$ using this model.

Then, if certain conditions are met, the following update rule for our proposed POCO algorithm is used.

\begin{definition}[Predictive update]
Let $\beta_t > 0 $ be an appropriately chosen step size. The predictive update is
\begin{equation}
\xx_{t+1} = \proj_{\XX} \left( \overline{\xx}_{t+1} - \beta_t \gt\right). \label{eq:_proj_update_p}
\end{equation}
\label{def:pred_up}
\end{definition}
\vspace{-0.75cm}

The predictive update is to be used directly after the OCO update and will lead to a strict improvement over the OCO update under certain conditions. The aforementioned conditions will be discussed in the next sections and depend on the properties of the loss function. If the conditions are not met, $\oxt$ is directly used. Let $\delta > 0$ be the desired improvement when using the predictive update. We define the counter $c_t$: 
\[
c_{t+1} = \begin{cases}
c_{t} + 1 & \text{if } \left\|\xxt - \oxt \right\| \geq \delta\\
c_{t} & \text{otherwise}
\end{cases}
\]
with $c_0 = 0$. The variable $c_t$ represents the number of predictive updates as described in Definition~\ref{def:pred_up}. Let $\nu = c_T/T$ be the ratio of rounds using the predictive update to the total number of rounds.

Depending on the loss function, any regret-bounded OCO update can be used in the POCO framework. Back to the \texttt{OGD} example, the predictive \texttt{OGD} uses the update~\eqref{eq:ogd_update} and if certain conditions are met, $\xx_{t+1}$ is provided by~\eqref{eq:_proj_update_p} and if not, $\xxt = \oxt$. We write $\xx_{t+1}(\beta_t) = \proj_{\XX} \left( \overline{\xx}_{t+1} - \beta_t \gt\right)$ as a function of the step size $\beta_t > 0$ and let $\dt = \xx_{t+1}(\beta_t) - \oxt$ be the descent direction.

Next, we provide sufficient conditions for the estimated gradient $\gt$ to be a feasible descent direction. 
Later, we consider the step size selection problem. Particularly, two cases are considered where (i) the step sizes are constant and chosen \emph{a priori} based on a property of the loss functions, or (ii) the step sizes are selected through the application of a backtracking line search that enforces a modified online version of the Armijo condition~\cite{bertsekas2015convex}.

The following lemma introduces a sufficient condition for the estimated gradient $\gt$ to be a descent direction of the OCO problem~\eqref{eq:oco_problem}.

\vspace{-.2cm}
\begin{lemma}[Estimated descent direction] 
The vector $-\gt$ provided by the $\epsilon$-forecaster is a descent direction for $f_{t+1}(\oxx_{t+1})$ if $\| \gt \| > \epsilon$.
\label{lem:descent}
\end{lemma}
\vspace{-.2cm}
The proof of Lemma~\ref{lem:descent} is presented in Appendix~\ref{app:lem_desc}. The next lemma is adapted from~\cite{bertsekas2015convex} and ensures that the predictive step follows a feasible descent direction.

\vspace{-.2cm}
\begin{lemma}[Feasible estimated descent direction]
For all $\beta_t > 0$ and $\oxx_{t+1} \in \XX$, if $\|\gt\| > \epsilon$ and $\xx_{t+1}(\beta_t) \neq \oxx_{t+1}$ , then $\xx_{t+1}(\beta_t) - \oxx_{t+1}$ is a feasible descent direction at $\oxt$ and
$
\gt \Tr \left( \xx_{t+1}(\beta_t) - \oxx_{t+1} \right) \leq -\frac{1}{\beta_t} \left\| \xxt(\beta_t) - \oxt \right\|^2.
$
\label{lem:feas_descent}
\end{lemma}
\vspace{-0.2cm}

Similarly, the proof of Lemma~\ref{lem:feas_descent} is given in Appendix~\ref{app:lem_feas_d}.

\section{POCO with fixed step size}
\label{sec:predetermined}

We now present a predictive update where step sizes $\beta_t$ are fixed and based on a propriety of the sequence of loss functions. We conclude this section by providing regret-bounded algorithms using these updates. In this section, we add the following assumption:

\vspace{-.2cm}
\begin{assumption} 
Let $L < \infty$. The loss function $f_t(\xx)$ has an $L$-uniformly Lipschitz-continuous gradient:
$
\left\|\nabla f_t(\xx) - \nabla f_t(\mathbf{y})\right\| \leq L \|\xx - \mathbf{y} \|
$
for all $t=1,2,\ldots, T$ and $\xx, \mathbf{y} \in \XX$.
\label{assumption:L_grad}
\end{assumption}
\vspace{-.2cm}

We propose a predictive update with fixed step size next. We state sufficient conditions that guarantee a strict improvement over an OCO update. These sufficient conditions can be checked at each round to determine if the estimated information is accurate enough, and therefore if the predictive update should be used in the current round.

\begin{lemma}[Predictive update with fixed step size]
Suppose that Assumption~\ref{assumption:L_grad} holds and $\| \gt \| > \epsilon$. If $\beta \leq \frac{1}{L}$ and $\| \dt \|= \|\proj_{\XX} \left( \overline{\xx}_{t+1} - \beta_t \gt\right) - \oxt \| \geq \frac{\epsilon}{L} + \sqrt{\frac{\epsilon^2}{L^2} + \frac{2\delta}{L}}$, then the predictive update~\eqref{eq:_proj_update_p} used by the $\epsilon$-forecaster strictly improves on the OCO update and the improvement is bounded below by $\delta > 0$.
\label{lem:proj_fixed_step_update}
\end{lemma}
\vspace{-.2cm}

The proof of Lemma~\ref{lem:proj_fixed_step_update} is provided in Appendix~\ref{app:lem_fixed}. We now present regret bounds for POCO algorithms. This algorithm uses the predictive update with fixed step size to improve the performance of OCO algorithms.

\begin{theorem}[POCO regret bound]
Consider an OCO algorithm with a sublinear regret upper bound. Suppose that the forecaster uses the predictive update~\eqref{eq:_proj_update_p} only at rounds $t$ when the estimated gradient $\gt$ and feasible descent direction $\dt = \proj_{\XX} \left( \overline{\xx}_{t+1} - \beta_t \gt\right) - \oxt$ satisfy the assumptions of Lemma~\ref{lem:proj_fixed_step_update}. If the ratio of rounds satisfying these assumptions is greater than $\nu$, then the regret of the POCO algorithm is bounded above by
\[
\emph{\reg}^d_T(POCO) \leq \emph{\reg}^d_T(OCO) - T \nu \delta.
\]
\label{thm:poco}
\end{theorem}
\vspace{-.8cm}

The proof of Theorem~\ref{thm:poco} is presented in Appendix~\ref{app:poco}. This theorem leads to the following corollary which provides a regret bound for the \texttt{OGD} with predictive updates (\texttt{POGD}).

\vspace{-.2cm}
\begin{corollary}[$O\left(\sqrt{T}\right)$ regret bound for \texttt{POGD}]
Suppose that the ratio $\nu$ of rounds that respects the assumptions of Lemma~\ref{lem:proj_fixed_step_update} is $\nu > \frac{1}{\sqrt{T}}$. Then the predictive \texttt{OGD} algorithm's regret is bounded above by
\begin{align*}
\emph{\reg}^d_T(\texttt{POGD}) &\leq \emph{\reg}^d_T(\texttt{OGD}) - \delta \sqrt{T},\\
&= \left(\frac{7X^2}{4} + \frac{G^2}{2} + X V_T - \delta \right) \sqrt{T},
\end{align*}
which is sublinear and tighter than the \texttt{OGD} regret bound.
\label{cor:pogd}
\end{corollary}
\vspace{-0.2cm}

The corollary follows from substituting $\reg^d_T(\texttt{OGD)}$ from~\cite{zinkevich2003online} and $\nu > 1/\sqrt{T}$ in Theorem~\ref{thm:poco}.

\section{POCO with backtracking line search}
\label{sec:backtracking}

In this section, we do not require Assumption~\ref{assumption:L_grad} to hold. We however use the following proposition:

\vspace{-0.2cm}
\begin{proposition}
The loss function $f_t(\xx)$ is $\Delta$-time-Lipschitz with $\Delta_t(\xx), \Delta < \infty$, that is:
\[
\left| f_t(\xx) - f_{t+\tau}(\xx) \right| \leq \Delta_t(\xx) |\tau| \leq \Delta |\tau|
\]
for all $\tau \in \left\{\left. i \in \mathbb{Z} \right| 0 \leq t+i \leq T \right\}$ at all $t=1,2,\ldots, T$, and all $\xx \in \XX$.
\label{assumption:time_Lip}
\end{proposition}
\vspace{-0.2cm}

Proposition~\ref{assumption:time_Lip} always holds because Assumption~\ref{assumption:bounded} implies that $\Delta = 2 B/\left|\tau\right|$ is sufficient.
Under Proposition~\ref{assumption:time_Lip}, we consider functions that are $(t,\xx)$-locally and globally Lipschitz in their time argument, respectively, for the intermediary bound ($\Delta_t(\xx)$) and the upper bound ($\Delta$). This can represent, for example, loss functions like squared tracking error functions, in which the time-varying targets are always contained in a closed set. 

In the case of the POCO with backtracking (POCOb), we re-express the update~\eqref{eq:_proj_update_p} given in Definition~\ref{def:pred_back}. The backtracking line search for predictive update is given in Algorithm~\ref{alg:proj_backtracking}. 

\vspace{-0.2cm}
\begin{definition}[POCOb update]
Let $\zeta$ be a positive scalar and $\beta^m$ be determined by a backtracking line search algorithm. The predictive update with backtracking line search is
\begin{equation}
\xxt = \oxt + \beta^m  \left( \proj_{\XX} \left( \overline{\xx}_{t+1} - \zeta \gt\right)  - \oxt \right)
\label{eq:update_backtracking}
\end{equation}
\label{def:pred_back}
\end{definition}
\vspace{-0.9cm}

\begin{algorithm}[tb]
\begin{algorithmic}[1]
\STATE \textbf{Parameters:} Given $\beta \in ( 0 , 1 )$ and $M \in \mathbb{N}$.
\STATE \textbf{Initialization:} Set $\zeta>0$.
\medskip

\STATE $\dt = \proj_\XX \left(\oxt - \zeta \gt \right) - \oxt$
\STATE $m=0$.
\WHILE{$f_t\left(\oxx_{t+1} + \beta^m \dt \right) > f_t(\oxt) + \beta^m \left(\gt \Tr \dt - \epsilon \|\dt\|\right) - 2 \Delta$ and $m \leq M$}
\STATE $m = m+ 1$.

\ENDWHILE

\medskip

\IF{$m > M$} 

\STATE $\beta = 0$.

\ENDIF

\end{algorithmic}
\caption{Backtracking algorithm for predictive gradient projection}
\label{alg:proj_backtracking}
\end{algorithm}

The next lemma shows that the backtracking line search-based predictive update improves on the OCO update. Our claim relies on the modified Armijo condition for gradient projection. This condition ensures a sufficient decrease in the objective when using an estimated gradient projection descent direction~\cite{wright1999numerical}. We adapt this condition to the estimated gradient and online setting. The modified Armijo condition for gradient projection~\cite{bertsekas2015convex} on $f_{t+1}$ and feasible descent direction $\dt=\proj_\XX \left(\oxt + \zeta \gt \right) - \oxt$ for some $\zeta>0$ with step size $\beta^m$ is given by:
\begin{align}
f_{t+1}\left(\oxx_{t+1} + \beta^m \dt \right) &\leq f_{t+1}(\oxx_{t+1}) \label{eq:proj_armijo}\\
&\quad+ \beta^m \nabla f_{t+1}(\oxx_{t+1})\Tr \dt. \nonumber
\end{align}

\begin{lemma}[Sufficient decrease of POCOb update]
Suppose $\| \gt \| > \epsilon$. If Algorithm~\ref{alg:proj_backtracking} terminates to a step size $\beta^m > 0$, then the predictive update with backtracking line search~\eqref{eq:update_backtracking} used by the $\epsilon$-forecaster satisfies the modified Armijo condition~\eqref{eq:proj_armijo}, and leads to a sufficient decrease in the loss function.
\label{lem:proj_backtracking}
\end{lemma}
\vspace{-0.2cm}
The proof of the lemma can be found in Appendix~\ref{app:lem_back}.

\vspace{-0.2cm}
\begin{remark}
Algorithm~\ref{alg:proj_backtracking} ensures that when $\beta \neq 0$, $\beta^m$ satisfied: 
\begin{align}
f_{t}\left(\oxx_{t+1} + \beta^m \dt \right) &\leq f_{t}(\oxx_{t+1}) + \beta^m \gt \Tr \dt \nonumber\\
&\quad- \beta^m \epsilon \|\dt\|  - 2\Delta \label{eq:proj_negative}
\end{align}
Every element of~\eqref{eq:proj_negative} is available at time $t$, which is not the case in~\eqref{eq:proj_armijo}. This allows us to use a backtracking line search algorithm to determine $\beta_t$ in an OCO setting. Algorithm~\ref{alg:proj_backtracking} also ensures that the step size is not too small (cf.~\cite[Section 3.1]{wright1999numerical}).
\end{remark}
\vspace{-0.2cm}

Note that there is an additional $\epsilon \| \dt \|$ term in the modified Armijo condition for estimated gradient projection. This is a consequence of not having access to the exact gradient of $f_t$. Hence, to ensure that the update is valid, the modified Armijo condition is augmented by a term proportional to the error of the estimated gradient. The second additional term, $2 \Delta$, is due to the time-varying setting of OCO.

We now discuss the existence of step sizes $\beta_t$ that satisfy~\eqref{eq:proj_negative} at round $t$. Before stating the main result, for a given $\oxt$ and $\gt$, define the set of step sizes that comply with line 5 in the line search algorithm, which is the modified Armijo condition for online settings~\eqref{eq:proj_negative}:
\begin{align*}
\mathcal{S} = \Big\{ \beta >0 \Big|& f_{t}\left(\oxx_{t+1} + \beta \dt \right) \leq f_{t}(\oxx_{t+1}) \\
&+ \beta \gt \Tr \dt - \beta \epsilon \|\dt\|  - 2\Delta \Big\}.
\end{align*}

\vspace{-0.8cm}
\begin{theorem}
Suppose $\dt = \proj_{\XX} \left( \overline{\xx}_{t+1} - \zeta \gt\right)  - \oxt \neq \mathbf{0}$ is a feasible descent direction and $f_t$ is bounded below for all $t$. Then there exists $\xx \in \XX$ such that $f_t(\oxt) - f_t(\xx) > 2 \Delta$ if and only if $\mathcal{S} \neq \emptyset$.
\label{thm:existence}
\end{theorem}
\vspace{-0.2cm}

The proof for Theorem~\ref{thm:existence} is presented in Appendix~\ref{app:existence}. We note that Theorem~\ref{thm:existence} does not guarantee that the backtracking algorithm, Algorithm~\ref{alg:proj_backtracking}, will find a non-zero step size. Other techniques like exact line searches, might be required to identify an adequate step size in some problem instances. Using Theorem~\ref{thm:existence}, we can provide a lower bound on the improvement of the predictive update with backtracking line search.

\vspace{-0.2cm}
\begin{corollary}[POCOb update improvement]
Suppose that the assumptions of Lemma~\ref{lem:proj_backtracking} hold and $\beta^m > 0$, then the predictive update with backtracking line search improves on the OCO update by a minimum of $2\Delta$.
\label{cor:backtracking_min}
\end{corollary}
\vspace{-0.2cm}

The proof of Corollary~\ref{cor:backtracking_min} is given in Appendix~\ref{app:pocob_min}. We now state a regret bound for the POCOb algorithm.

\vspace{-0.2cm}
\begin{theorem}[POCOb regret bound]
Consider an OCO algorithm with bounded regret. Suppose that the assumptions of Lemma~\ref{lem:proj_backtracking} are met. If the ratio of rounds with $\beta > 0$ and satisfying these assumptions to $T$ is greater than $\nu$, then the regret of the POCO algorithm with backtracking used by the $\epsilon$-forecaster is bounded  above by
\begin{equation}
\emph{\reg}^d_T(POCOb) \leq \emph{\reg}^d_T(OCO) - 2 T \nu \Delta
\label{eq:regret_pocob}
\end{equation}
and thus outperforms the OCO algorithm.
\label{thm:pocob}
\end{theorem}
\vspace{-0.2cm}

The proof of Theorem~\ref{thm:pocob} is presented in Appendix~\ref{app:pocob}.

\begin{remark}
Note that if the locally Lipschitz statement of Proposition~\ref{assumption:time_Lip} is used, then $2 \Delta$ is replaced by $\Delta_{t,1}(\oxt + \beta^m \dt) + \Delta_{t,1}(\oxt)$ in the modified Armijo condition for online settings~\eqref{eq:proj_negative}, and the bound~\eqref{eq:regret_pocob} can be recomputed accordingly.
\end{remark}

\vspace{-0.25cm}

\section{Example}
\vspace{-0.25cm}

\label{sec:example}
In this section, we apply POCO algorithms to demand response (DR) in power systems~\cite{callaway2011achieving,palensky2011demand}, specifically regulation and curtailment. At each time step, a DR aggregator sends instructions to their loads to follow a regulation signal, e.g., a power imbalance due to a sudden change in renewable power generation~\cite{callaway2009tapping,taylor2016power}. Each load responds to the signal by adjusting its power consumption. The power consumption is constrained by a storage capacity, which could represent physical storage like a battery or the load's limits, e.g., thermal constraints. The regulation signal is unknown at the time the DR instructions are sent. This can be due, for example, to a drop in renewable power generation which is only assessed after the generator has committed to some amount of power. The objective of the DR aggregator is, therefore, to predict the DR dispatch at each time instance. This problem can be formulated as POCO, in which an estimate of the regulation signal is available to the load aggregator.

We consider $N$ loads. Let $\xx_t \in \RR^N$ denote the decision variable at round $t$. The variable $\xx_t$ represents the instructions sent to the loads. Let $r_t \in \RR$ be the regulation signal at time $t$. Let $\overline{\xx}, \underline{\xx} \in \RR^N$ be the maximum and minimum power that can be consumed or delivered for all loads. Define the decision set $\XX = \left\{\left. \xx \in \RR^N \right| \underline{\xx} \leq \xx \leq \overline{\xx} \right\}$. We let $\ss_{t} \in \RR^N$ denote the state of charge vectors of the loads at time $t$ and $\cc \in \RR^N$ the vector of load energy capacities. The state of charge of a load $i$ at time $t$ is $s_t(i) = s_0(i) + \sum_{n=1}^t x_n(i)$. In the current case, we assume that there is no leakage nor energy losses. The OCO problem takes the following form:
\begin{equation}
\min_{\xx_t \in \XX} \left( r_t - \mathbf{1}\Tr \xx_t \right)^2 + \sigma \left\| \ss_{t-1} + \xx_t - \frac{\mathbf{c}}{2} \right\|^2.
\label{eq:oco_reg}
\end{equation}
The loss function has two terms: (i) a regulation term where the aggregated loads are dispatched to follow a regulation signal $r_t$ and (ii) a state of charge objective added to keep the loads near half their energy capacity. The loss function given in~\eqref{eq:oco_reg} is $\sigma$-strongly convex. For this reason we use the \texttt{OGD} for strongly convex functions ($\sigma$\texttt{OGD}) proposed in~\cite{mokhtari2016online}, which offers tighter regret bound than the standard \texttt{OGD}. The following corollary gives an upper bound on the regret of predictive \texttt{OGD} for strongly convex function ($\sigma$\texttt{POGD}).

\vspace{-0.2cm}
\begin{corollary}[\texttt{POGD} for strongly convex functions]
Suppose $f_t$ is $\sigma$-strongly convex  and satisfies Assumption~\ref{assumption:L_grad} for all $t$.
Consider the $\sigma$\texttt{OGD} update
\begin{align*}
\overline{\xx}_{t+1} &= \xx_{t} + \eta  \left( \proj_{\XX} \left( \xx_{t} - \frac{1}{\gamma} \nabla f_{t}(\xx_{t})\right)  - \xx_{t} \right)
\end{align*}
where $\eta \in (0,1]$ and $0 < \gamma \leq L$. Then, the $\sigma$\texttt{POGD} with fixed step size, given that the assumptions of Lemma~\ref{lem:proj_fixed_step_update} hold for a ratio of the total rounds greater than $\nu$, has a regret bounded above by
\begin{align*}
  \emph{\reg}^d_T(\sigma \texttt{POGD}) &\leq \emph{\reg}^d_T(\sigma\texttt{OGD}) - T \nu \delta,\\
  &\leq O\left( V_T + 1 \right) - T \nu \delta.
\end{align*}
\label{cor:strongly}
\end{corollary}
\vspace{-1cm}

The result follows from the proof of Theorem~\ref{thm:poco} and the $\sigma$\texttt{OGD} regret bound from~\cite{mokhtari2016online}. We now present simulation results. All optimizations are solved using \texttt{CVXPY}~\cite{cvxpy} and the \texttt{ECOS}~\cite{ecos} solver. 

\vspace{-0.75cm}
\paragraph*{Fixed step size example.}
\label{ssec_fixed}

The load and algorithm parameters for this example are gathered in Table~\ref{tab:poco_para}. The initial state of charge of each load is set to half its capacity. The regulation signal is $r_t = 0.2 \sin \left(\frac{2\pi}{T} t \right) + w_t$. The parameter $w_t \sim \mathrm{N}(0,0.01)$ is a Gaussian noise used to model sudden changes. We assume that the aggregator has access to estimated gradient for different level of accuracy $\epsilon$. This represents, for example when $\epsilon = 0.01$, a relative error of at least $4\%$ of the actual gradient norm. The parameter $\sigma$ is set to achieve adequate regulation performance without deviating too much from each load's desired state of charge.

\begin{table}[tb]
  \caption{Parameters for POCO simulations}
  \label{tab:poco_para}
  \renewcommand{\arraystretch}{1.3}
  \centering

  \vspace{0.1cm}
  \begin{tabular}{ccc}
  \hline

  \hline
  \textbf{Parameter} & \textbf{Value} & \textbf{Unit}\\
  \hline
  $N$   & $25$ & loads\\
  $h$   & $30$ & seconds\\
  $\overline{\xx}/h$ & $\mathrm{Uniform}[1,3]$ & kW\\
  $\underline{\xx}/h$ & $-\overline{\xx}$ & kW\\
  $\mathbf{c}$ & $\mathrm{Uniform}[10,15]$ & kWh\\
  $\epsilon$ & $0.1$, $0.05$ \& $0.01$ & --- \\
  $\delta$ & $10^{-6}$ & --- \\
  $\sigma$ & $0.005$ & --- \\
  $\eta$ & $1$ & --- \\
  $\gamma$ & $L$ & --- \\
  $\beta$ & $1/L$ & --- \\
  \hline

  \hline
  \end{tabular}
\end{table}

\begin{figure}[ht]
\vspace{-0.2cm}
  \centering
  \includegraphics[width=1\columnwidth]{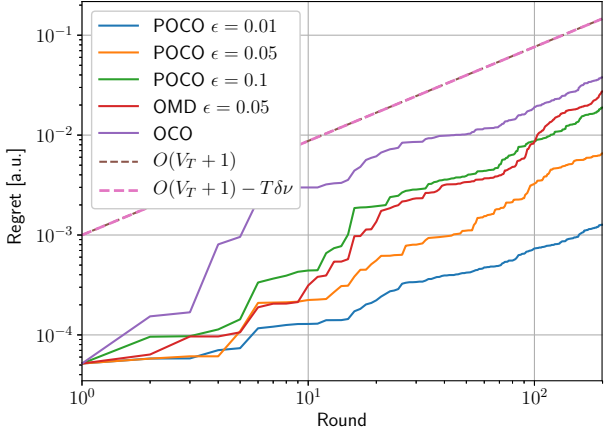}
  \vspace{-0.75cm}
  \caption{Regret comparison between POCO with fixed step size, OCO, and \texttt{OMD} (log scale, note: $O\left(V_T +1\right)$ and $O\left(V_T +1\right) - T\delta\nu$ are superimposed when shown at this scale)}
  \label{fig:regret_log}
\end{figure}

We now present the performance of our POCO algorithm with a fixed step size. We implement the \texttt{OMD} from~\cite{rakhlin2013optimization} for comparison. This algorithm uses $\gt$ without validating the estimated information. Figure~\ref{fig:regret_log} shows an instance of the experimental regret for the POCO with three different values of $\epsilon$, the conventional OCO algorithm, their respective regret bounds and \texttt{OMD}'s regret. POCO outperforms its bound, the OCO and the \texttt{OMD} algorithms. We remark that as expected the number of predictive updates increases with the accuracy of the estimated gradient, the performance of the POCO algorithm also improves.
\vspace{-0.75cm}

\paragraph*{Backtracking line search example.}
\label{sssec_cur}
We now present an example of POCO with backtracking. We consider a curtailment scenario. We let $p_t$ be the total power to be curtailed by the loads at time $t$ for $t = 1, 2, \ldots, T$. When a contingency occurs in the network, flexible loads are called to curtail their power consumption, e.g., by temporarily shutting down their HVAC system. Contrary to the regulation case, the loads are not contracted to follow a setpoint and no penalties are assessed on loads curtailing more than asked. Similar to the regulation setting, the curtailment signal is unknown until immediately after the current round. This setting can be modeled as POCO where an estimated curtailment signal is available to the aggregator at each round.
We use the same notation as the previous examples. Let $[\cdot]^+ = \max\{0, \cdot\}$. This curtailment scenario is modeled by loss function given below:
\begin{equation}
f_t(\xx_t) = 
\left(\left[p_t - \mathbf{1}\Tr \xx_t \right]^+\right)^2 + \sigma \left\| \alpha \ss_{t-1} + \xx_t - \frac{\mathbf{c}}{2}  \right\|^2
\label{eq:f_not_Lip}
\end{equation}
where we have added a recovery coefficient to the state of charge objective term used previously. This coefficient models the usual evolution of the load (e.g., ambient temperature heating for a thermostatic load). We let $\alpha = 1.001$. This is equivalent to a recovery coefficient of $1.13$ per hour. The function $f_t$ given in~\eqref{eq:f_not_Lip} is not gradient Lipschitz and Assumption~\ref{assumption:L_grad} does not hold. We model the curtailment signal to be quickly increasing at first and then slowly plateauing to represent new level of available generation. This event is assumed to be limited in time, after which the network goes back to its normal state and no curtailment is then required. We let $p_t = 0.04 t^{0.3} + w_t$ where $w_t \sim \mathrm{N}(0,0.01)$ for $t=1,2, \ldots, T/4 $ and then $p_t = 0.04 \left(T/4\right)^{0.3} + w_t'$ where $w_t' \sim \mathrm{N}(0,0.001)$ for $t=T/4, T/4+1, \ldots, T $. The noise variance is equivalent to approximatively $10\%$ of curtailment signal's value at first and then about $1\%$.

\begin{table}[tb]
  \caption{Different parameters for POCOb simulations}
  \label{tab:pocob_para}
  \renewcommand{\arraystretch}{1.3}
  \centering

  \begin{tabular}{ccc}
  \hline

  \hline
  \textbf{Parameter} & \textbf{Value}\\
  \hline
  $\alpha$   & $1.001$\\
  $\epsilon$ & $0.1$, $0.01$ \& $0.001$ \\
  $M$ & $100$ \\
  $\sigma$ & $5 \times 10^{-5}$ \\
  $\zeta$ & $0.5$ \\
  $\beta$ & $0.9$ \\
  $\eta$ & $1/10\sqrt{T}$ \\
  \hline

  \hline
  \end{tabular}
\end{table}

\begin{figure}[tb]
  \vspace{-0.21cm}
  \centering
  \includegraphics[width=1\columnwidth]{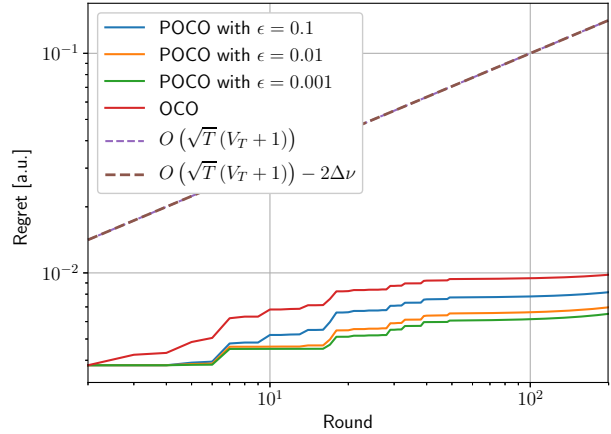}
  \vspace{-0.75cm}
  \caption{Regret comparison between POCOb, OCO, and \texttt{OMD} (log scale, note: $O\left(\sqrt{T}\left(V_T +1\right)\right)$ and $O\left(\sqrt{T}\left(V_T +1\right)\right) - T\delta\nu$ are superimposed when shown at this scale)}
  \label{fig:regret_log_cur}
  \vspace{-0.05cm}
\end{figure}

We use the same parameters as in the previous section, except for those in Table~\ref{tab:pocob_para}. The POCOb experimental regret shown in Figure~\ref{fig:regret_log_cur} is sublinear in the numbers of rounds and outperforms the OCO's regret. While the performance is not as strong as POCO with fixed step size, this algorithm can be applied to a broader family of functions because it does not require the loss function to be gradient Lipschitz continuous. The difference in performances between the two POCO updates is explained by the fact that the sufficient conditions for the POCOb are rarely satisfied in this simulation. Improved estimation accuracy, $\epsilon$, and a loss function with a lower maximum temporal change, $\Delta$, could increase the number of times the backtracking line search is used. Nevertheless, the POCOb achieves a regret reduction of 29\% when $\epsilon = 1\%$ over a standard OCO algorithm. Lastly, we note that POCOb performs better for larger variations in $p_t$ and smaller values of $\epsilon$.

\section{Conclusion}
In this work, we have presented the predictive online convex optimization framework. In POCO, a second update is used after the OCO update to improve performance using an estimated gradient. We have presented two versions of the predictive update that can be used under different assumptions. We have shown a regret upper bound for all of our POCO algorithms. We have applied POCO to demand response in electric power systems and found that they outperform conventional OCO using commonly available forecast information. In the case of fixed step size update, we observed an improvement of $95\%$ in the final regret and of $29\%$ in the backtracking case when having access to a $(\epsilon=0.01)$-forecaster.

\vspace{-0.4cm}
\begin{ack} 
\vspace{-0.35cm} 

A.L.-L. thanks P. Mancarella for his support and for co-hosting him at The University of Melbourne.
\end{ack}
\vspace{-0.1cm}

\appendix
\vspace{-0.3cm}
\section{Proof of Lemma~\ref{lem:descent}}
\label{app:lem_desc}
\vspace{-0.4cm}

Define $\ee_t \in \RR^n$ as $\ee_t = \gt - \nabla f_{t+1}(\overline{\xx}_{t+1})$ where $\|\ee_t \| \leq \epsilon$ by the definition of the $\epsilon$-forecaster. A vector $\mathbf{v}$ is a descent direction if $\mathbf{v}^\top \nabla f_{t+1}(\overline{\xx}_{t+1}) < 0$. Thus, we have
\begin{align*}
0 &< \gt\Tr \nabla f_{t+1}(\overline{\xx}_{t+1}) \\
&= \gt\Tr (\gt - \ee_t)
\end{align*}
Equivalently, we have $\gt\Tr\gt > \gt\Tr \ee_t$. Taking the norm of both sides and dividing by the norm of $\gt$ gives
\begin{equation}
\|\gt\| > \|\ee_t\| \cos \theta_{\gt,\ee_t},\label{eq:should_hold}
\end{equation}
where $\theta_{\gt,\ee}$ is the angle between $\gt$ and $\ee_t$.
By assumption, $\|\gt\| > \epsilon$ and $\| \ee_t \| \leq \epsilon$. Therefore~\eqref{eq:should_hold} always holds and we have proved the lemma.  \QED

\vspace{-0.3cm}
\section{Proof of Lemma~\ref{lem:feas_descent}}
\label{app:lem_feas_d}
\vspace{-0.4cm}

The identity follows from~\cite[Proposition 6.1.1]{bertsekas2015convex} with $\gt$ instead of the gradient of the loss function. It then follows from Lemma~\ref{lem:descent} that $-\beta_t\gt$ with $\beta_t >0$ is a descent direction at $\oxt$. Thus, $\dt = \xx_{t+1}(\beta_t) - \oxx_{t+1}$ is a feasible descent direction because $\dt \in \XX$ and $\gt \Tr \dt < 0$ for all $t$ and $\oxt \in \XX$. \QED

\vspace{-0.3cm}
\section{Proof of Lemma~\ref{lem:proj_fixed_step_update}}
\label{app:lem_fixed}
\vspace{-0.4cm}
By Assumption~\ref{assumption:L_grad}, $f_{t+1}$ has an $L$-Lipschitz gradient. We use the following inequality from ~\cite[Theorem 2.1.5]{nesterov1998introductory}
\vspace{-0.2cm}
\begin{equation}
\begin{aligned}
f_{t+1}(\mathbf{y}) &\leq f_{t+1}(\xx) + \nabla f_{t+1} (\xx)\Tr (\mathbf{y} - \xx)\\
&\quad+ \frac{L}{2} \left\| \xx - \mathbf{y} \right\|^2
\end{aligned}
\label{eq:lip_grad}\vspace{-0.3cm}
\end{equation}
for all $\xx, \mathbf{y} \in \XX$. We substitute $\mathbf{y} = \xxt(\beta)$ and $\xx = \oxx_{t+1}$ into~\eqref{eq:lip_grad} to obtain
\vspace{-0.2cm}
\begin{align*}
f_{t+1}(\xxt(\beta)) &\leq f_{t+1}(\oxt) + \frac{L}{2} \left\|\xxt(\beta)- \oxt \right\|^2\\
&\quad + \nabla f_{t+1} (\oxt)\Tr (\xxt(\beta) - \oxt). \vspace{-0.3cm}
\end{align*}
For the reminder of the proof, we use $\dt = \xxt(\beta)- \oxt$ to simplify the notation. We rewrite the gradient in term of the estimated gradient, which yields
\begin{align}
f_{t+1}(\xxt(\beta)) &\leq f_{t+1}(\oxt) + \gt\Tr \dt-\ee_t\Tr \dt \nonumber \\
&\quad+ \frac{L}{2} \left\|\dt \right\|^2. \label{eq:l_sub} \vspace{-0.5cm}
\end{align}
By assumption, $\xxt(\beta) \neq \oxt$, which ensures that Lemma~\ref{lem:feas_descent} holds. We use Lemma~\ref{lem:feas_descent} to upper bound the second term of the right-hand side of~\eqref{eq:l_sub}. We then have:
$f_{t+1}(\xxt(\beta))\leq f_{t+1}(\oxt) -\left(\frac{1}{\beta}-\frac{L}{2}\right) \left\| \dt \right\|^2 
+ \epsilon \|\dt\|$.
Therefore, the predictive update with fixed step size will improve on the OCO update by a minimum of $\delta > 0$ if the following condition is satisfied:
\begin{equation}
\frac{1}{\beta}\left\| \dt \right\|^2 -\frac{L}{2} \left\| \dt \right\|^2 - \epsilon \|\dt\| \geq \delta.
\label{eq:proj_star}
\end{equation}
Assuming $0 < \beta \leq \frac{1}{L}$, then $\frac{1}{\beta} \geq L$, and if $\frac{L}{2}\left\| \dt \right\|^2 - \epsilon \|\dt\| \geq \delta$, then~\eqref{eq:proj_star} also holds for any $\beta \in ] 0 , \frac{1}{L}]$. Solving for the norm of the feasible descent direction $\| \dt \|$, we have
\begin{equation}
\| \dt \| = \|\xxt(\beta) - \oxt\| \geq \frac{\epsilon}{L} + \sqrt{\frac{\epsilon^2}{L^2} + \frac{2\delta}{L}}. \label{eq:dir_cond}
\end{equation}
Thus, by setting $0 < \beta \leq \frac{1}{L}$ and satisfying~\eqref{eq:dir_cond}, we obtain $f_{t+1}(\xxt(\beta)) \leq f_{t+1}(\oxt) -\delta, \label{eq:fixed_conclu}$ where $\delta > 0$. This implies that the predictive update strictly improves over the OCO update when the feasible descent direction satisfies the condition~\eqref{eq:dir_cond}. The improvement is bounded below by $\delta$. \QED

\vspace{-0.3cm}
\section{Proof of Theorem~\ref{thm:poco}}
\label{app:poco}
\vspace{-0.4cm}

Let $\hat{\xx}_t$ denote the decision variable with $\beta_t = 0$ for all $t$. In other words, $\hat{\xx}_t$ represents the decision variable computed without the predictive algorithm. Denote the set of assumptions of Lemma~\ref{lem:proj_fixed_step_update} at round $t$ by $\mathcal{A}_t$. Let $\mathbb{I}_{\mathcal{A}_t}$ be the indicator function where $\mathbb{I}_{\mathcal{A}_t} =1$ if the assumptions are satisfied and $0$ otherwise. Observe that the improvement, $i_t$, is given by
\begin{equation}
i_t \mathbb{I}_{\mathcal{A}_t} = f_t(\hat{\xx}_t) - f_t(\xx_t),
\label{eq:improv}
\end{equation}
where $i_t$ is the improvement when $\mathbb{I}_{\mathcal{A}_t}=1$. The regret of the POCO algorithm is
\begin{equation}
\reg^d_T(POCO) = \sum_{t=1}^T f_t(\xx_t) - f_t(\xx_t^\ast).
\label{eq:regret_poco}
\end{equation}
Using~\eqref{eq:improv}, we re-express $f_t(\xx_t)$ in~\eqref{eq:regret_poco}:
\begin{align}
\reg^d_T(POCO) &= \sum_{t=1}^T f_t(\hat{\xx}_t) - f_t(\xx_t^\ast) - \sum_{t=1}^T i_t \mathbb{I}_{\mathcal{A}_t} \nonumber\\
&= \reg^d_T(OCO) - \sum_{t=1}^T i_t \mathbb{I}_{\mathcal{A}_t} \label{eq:upper_temp}
\end{align}
By Lemma~\ref{lem:proj_fixed_step_update}, the improvement $i_t$ is bounded below by $\delta$. We rewrite~\eqref{eq:upper_temp} as
$\reg^d_T(POCO) \leq \reg^d_T(OCO) - \sum_{t=1}^T \delta \mathbb{I}_{\mathcal{A}_t}.$
A minimum of $T\nu$ rounds satisfy $\mathcal{A}_t$ and hence $\reg^d_T(POCO) \leq \reg^d_T(OCO) - T\nu\delta$. \QED

\vspace{-0.3cm}
\section{Proof of Lemma~\ref{lem:proj_backtracking}}
\label{app:lem_back}
\vspace{-0.4cm}

We show that for some step size $\beta^m$, the estimated gradient descent projection leads to a sufficient decrease thus outperforming the OCO update. We show that if $\beta^m$ satisfies~\eqref{eq:proj_negative}, then it also satisfies satisfies~\eqref{eq:proj_armijo}, ensuring a sufficient decrease in the objective function. Note that~\eqref{eq:proj_negative} is the condition under which the backtracking algorithm, Algorithm~\ref{alg:proj_backtracking}, is used.
We can see from the left-hand side of the condition~\eqref{eq:proj_negative} that the update improves over the OCO update because the three last terms are bounded above by $0$, i.e., $\gt \Tr \dt \leq -\frac{1}{\zeta} \| \dt\|^2 < 0$ by Lemma~\ref{lem:feas_descent}. Thus  all three terms are less or equal to zero. By assumption, $\beta>0$, and these terms are also bounded away from zero since $\xxt \neq \oxt$.

We start from~\eqref{eq:proj_negative} and shows it implies~\eqref{eq:proj_armijo}. By assumption, $\| \ee_t \| \leq \epsilon$ for all $t$ and hence~\eqref{eq:proj_negative} implies,
\begin{align}
f_{t}\left(\oxx_{t+1} + \beta^m \dt \right) &\leq f_{t}(\oxx_{t+1}) + \beta^m \left(\gt - \ee_t \right)\Tr \dt \nonumber\\
&\quad- 2\Delta. \nonumber
\end{align}
Rearranging the terms, we have
\begin{align}
f_{t}\left(\oxx_{t+1} + \beta^m \dt \right) + \Delta &\leq f_{t}(\oxx_{t+1}) - \Delta \label{eq:proj_starting} \\
&\quad+ \beta^m \nabla f_{t+1}(\oxx_{t+1})\Tr \dt. \nonumber
\end{align}
By assumption, $f_{t+1}$ is time-Lipschitz with constant $\Delta > 0$ for all $\xx \in \XX$ and all $t$. We can therefore bound below and above respectively the left-hand and right-hand side of~\eqref{eq:proj_starting}. This leads to
\begin{align*}
f_{t+1}\left(\oxx_{t+1} + \beta^m \dt \right) &\leq f_{t+1}(\oxx_{t+1}) \\
&\quad+ \beta^m \nabla f_{t+1}(\oxx_{t+1})\Tr \dt,
\end{align*}
the modified Armijo condition~\eqref{eq:proj_armijo}. \QED

\vspace{-.3cm}
\section{Proof of Theorem~\ref{thm:existence}}
\label{app:existence}
\vspace{-.4cm}
Assume $f_t(\oxt) - f_t(\xx) > 2 \Delta$. This assumption implies that $f_{t+1}(\oxt) - f_{t+1}(\xx) > 0$ by Proposition~\ref{assumption:time_Lip}.
Thus, $\oxt$ is not the minimum point of $f_{t+1}$. It follows that $\nabla f_{t+1} (\oxt) \neq 0$. By assumption, $\dt \neq \mathbf{0}$ is a feasible descent direction and we have
\begin{equation}
\nabla f_{t+1} (\oxt) \Tr \dt < 0.
\label{eq:descent_dir_def}
\end{equation}
Let $a \in (0, 1)$. Subtracting $\nabla f_t\left(\oxt + a \beta \dt\right)\Tr \dt$ on both side of~\eqref{eq:descent_dir_def} we obtain,
\begin{align}
( \nabla f_{t+1} (\oxt) - \nabla & f_t(\oxt + a \beta \dt) )\Tr \dt < \label{eq:subtracking}\\
 &-\nabla f_t\left(\oxt + a \beta \dt\right)\Tr \dt. \nonumber
\end{align}
If the following condition holds, then~\eqref{eq:subtracking} also holds:
\begin{equation}
\begin{aligned}
\| \nabla f_{t+1} (\oxt) &- \nabla  f_t(\oxt + a \beta \dt)\| \|\dt\|< 
\\&-\nabla f_t\left(\oxt + a \beta \dt\right)\Tr \dt.
\end{aligned}
\label{eq:upp_norm}
\end{equation}
Under Assumption~\ref{assumption:grad_B}, for all $\xx, \mathbf{z} \in \XX$ we have
\begin{align*}
\| \nabla f_{t+1} (\xx) - \nabla f_{t}(\mathbf{z})\| &\leq  \| \nabla f_{t+1} (\xx)\| + \|\nabla f_t(\mathbf{z})\| \leq 2G
\end{align*}
and by Assumption~\ref{assumption:compact}, we have $\|\dt\| \leq D$. Then, if 
\begin{equation}
2G D < -\nabla f_t\left(\oxt + a \beta^m \dt\right)\Tr \dt
\label{eq:upp_G}
\end{equation}
holds, so does~\eqref{eq:upp_norm}. We rewrite~\eqref{eq:upp_G} as
\begin{equation}
\nabla f_t\left(\oxt + a \beta \dt\right)\Tr \dt < - 2 G D
\label{eq:bound_grad}
\end{equation}

\vspace{-0.25cm}
Recalling Taylor's Theorem~\cite[Theorem 2.1]{wright1999numerical}:
$$f_t(\mathbf{y} + \mathbf{p}) = f_t(\mathbf{y}) + \nabla f_t (\mathbf{y} + a \mathbf{p}) \Tr \mathbf{p},$$
where $\mathbf{y}, \xx \in \XX$, $\mathbf{p}\in \RR^n$ and for some $a \in ( 0 , 1 )$. We let $\mathbf{y} = \oxt$ and $\mathbf{p} = \beta^m \dt$. We have,
\begin{align}
f_t(\oxt + \beta \dt) &= f_t(\oxt) \label{eq:appli_Taylor}\\
&\,+ \beta \nabla f_t (\oxt + a \beta \dt) \Tr \dt. \nonumber
\end{align}
We bound above the the last term of~\eqref{eq:appli_Taylor} using~\eqref{eq:bound_grad} and obtain
\begin{equation}
f_t(\oxt + \beta \dt) <  f_t(\oxt) - 2 \beta G D \label{eq:with_beta}
\end{equation}
By setting $\beta \leq \frac{\Delta}{G D}$ in~\eqref{eq:with_beta}, we have $f_t(\oxt + \beta \dt) <  f_t(\oxt) - 2 \Delta$. This shows that there always exists at least one point which satisfies the assumption on the existence of $\xx \in \XX$ such that $f_t(\oxt) - f_t(\xx) > 2 \Delta$ that is along the feasible descent direction $\dt$ from $\oxt$.

\vspace{-0.25cm}
Next, adapting the proof of~\cite[Lemma 3.1]{wright1999numerical} for the modified Armijo condition for online settings~\eqref{eq:proj_negative}, it follows that there exists $\overline{\beta} \leq \frac{\Delta}{GD}$ such that 
\begin{align*}
f_{t}\left(\oxx_{t+1} + \overline{\beta} \dt \right) &\leq f_{t}(\oxx_{t+1}) + \overline{\beta} \gt \Tr \dt - \overline{\beta} \epsilon \|\dt\|\\
&\quad  - 2\Delta.
\end{align*}
The set $\mathcal{S}$ is therefore non-empty if there exists $\xx \in \XX$ such that $f_t(\oxt) - f_t(\xx) > 2 \Delta$.
\vspace{-0.25cm}

We now show the converse. Assuming $\mathcal{S} \neq \emptyset$, then there exists $ \underline{\beta} \in \mathcal{S}$ and
\begin{align}
f_{t}\left(\oxx_{t+1} + \underline{\beta} \dt \right) &< f_{t}(\oxx_{t+1})  - 2\Delta ~\label{eq:beta_1}
\end{align}
holds since $\gt \Tr \dt < 0$ by Lemma~\ref{lem:feas_descent} and $\epsilon >0$. Thus,~\eqref{eq:beta_1} implies that there exists $\xx \in \XX$ such that $f_t(\oxt) - f_t(\xx) > 2 \Delta$ and one of such point is $\xx = \oxx_{t+1} + \beta_2 \dt$. This completes the proof. \QED

\vspace{-.3cm}
\section{Proof of Corollary~\ref{cor:backtracking_min}}
\label{app:pocob_min}
\vspace{-.4cm}

Since $\beta > 0$, then $\mathcal{S} \neq \emptyset$. By the converse of Theorem~\ref{thm:existence}, we have $f_{t}(\oxx_{t+1})  - f_t(\xx) > 2 \Delta,$ where $\xx = \oxt + \beta \dt$, the decision played by the predictive update~\eqref{eq:update_backtracking}. The predictive update hence improves on the OCO update by at least $2 \Delta$. \QED

\vspace{-.3cm}
\section{Proof of Theorem~\ref{thm:pocob}}
\label{app:pocob}
\vspace{-.4cm}

Let $\mathbb{I}_{\mathcal{A}'_t}$ be the indicator function where $\mathbb{I}_{\mathcal{A}'_t} =1$ if at round $t$, $\beta > 0$ and $\| \gt \| > \epsilon$ or $0$ otherwise. Using the same approach as in Theorem~\ref{thm:poco}'s proof with Corollary~\ref{cor:backtracking_min}, we obtain the regret bound. The last term of~\eqref{eq:regret_pocob} is strictly positive and thus the POCOb regret is always bounded above by the OCO algorithm regret. \QED


\vspace{-.5cm}

\begin{thebibliography}{10}
\vspace{-.5cm}

\bibitem{bertsekas1997nonlinear}
D.~P. Bertsekas.
\newblock Nonlinear programming.
\newblock {\em Journal of the Operational Research Society}, 48(3):334--334,
  1997.

\bibitem{bertsekas2015convex}
D.~P. Bertsekas.
\newblock {\em Convex optimization algorithms}.
\newblock Athena Scientific Belmont, 2015.

\bibitem{borrelli2017predictive}
F.~Borrelli, A.~Bemporad, and M.~Morari.
\newblock {\em Predictive control for linear and hybrid systems}.
\newblock Cambridge University Press, 2017.

\bibitem{callaway2009tapping}
D.~S Callaway.
\newblock Tapping the energy storage potential in electric loads to deliver
  load following and regulation, with application to wind energy.
\newblock {\em Energy Conversion and Management}, 50(5):1389--1400, 2009.

\bibitem{callaway2011achieving}
D.S. Callaway and I.~A. Hiskens.
\newblock Achieving controllability of electric loads.
\newblock {\em Proceedings of the IEEE}, 99(1):184--199, 2011.

\bibitem{cao2018virtual}
X.~Cao, J.~Zhang, and H.~V. Poor.
\newblock A virtual-queue based algorithm for constrained online convex
  optimization with applications to data center resource allocation.
\newblock {\em IEEE Journal of Selected Topics in Signal Processing}, 2018.

\bibitem{chen2018bandit}
T.~Chen and G.~B. Giannakis.
\newblock Bandit convex optimization for scalable and dynamic iot management.
\newblock {\em IEEE Internet of Things Journal}, 2018.

\bibitem{chen2017online}
T.~Chen, Q.~Ling, and G.~B. Giannakis.
\newblock An online convex optimization approach to proactive network resource
  allocation.
\newblock {\em IEEE Transactions on Signal Processing}, 65(24):6350--6364,
  2017.

\bibitem{chiang2012online}
C.-K. Chiang, T.~Yang, C.-J. Lee, M.~Mahdavi, C.-J. Lu, R.~Jin, and S.~Zhu.
\newblock Online optimization with gradual variations.
\newblock In {\em Conference on Learning Theory}, pages 6--1, 2012.

\bibitem{dekel2017online}
O.~Dekel, A.~Flajolet, N.~Haghtalab, and P.~Jaillet.
\newblock Online learning with a hint.
\newblock In {\em Advances in Neural Information Processing Systems}, pages
  5305--5314, 2017.

\bibitem{cvxpy}
S.~Diamond and S.~Boyd.
\newblock {CVXPY}: A {P}ython-embedded modeling language for convex
  optimization.
\newblock {\em Journal of Machine Learning Research}, 17(83):1--5, 2016.

\bibitem{ecos}
A.~Domahidi, E.~Chu, and S.~Boyd.
\newblock {ECOS}: {A}n {SOCP} solver for embedded systems.
\newblock In {\em European Control Conference (ECC)}, pages 3071--3076, 2013.

\bibitem{duchi2010composite}
J.~C. Duchi, S.~Shalev-Shwartz, Y.~Singer, and A.~Tewari.
\newblock Composite objective mirror descent.
\newblock In {\em COLT}, pages 14--26, 2010.

\bibitem{garcia1989model}
C.~E. Garcia, D.~M. Prett, and M.~Morari.
\newblock Model predictive control: theory and practice—a survey.
\newblock {\em Automatica}, 25(3):335--348, 1989.

\bibitem{hazan2016introduction}
E.~Hazan.
\newblock Introduction to online convex optimization.
\newblock {\em Foundations and Trends{\textregistered} in Optimization},
  2(3-4):157--325, 2016.

\bibitem{hazan2007logarithmic}
E.~Hazan, A.~Agarwal, and S.~Kale.
\newblock Logarithmic regret algorithms for online convex optimization.
\newblock {\em Machine Learning}, 69(2-3):169--192, 2007.

\bibitem{hazan2007online}
E.~Hazan and N.~Megiddo.
\newblock Online learning with prior knowledge.
\newblock In {\em International Conference on Computational Learning Theory},
  pages 499--513. Springer, 2007.

\bibitem{ho2016accelerating}
N.~Ho-Nguyen and F.~K{\i}l{\i}n{\c{c}}-Karzan.
\newblock \phantom{c}{A}ccelerating optimization under uncertainty via online
  convex optimization.
\newblock Technical report, August, 2016.

\bibitem{jadbabaie2015online}
A.~Jadbabaie, A.~Rakhlin, S.~Shahrampour, and K.~Sridharan.
\newblock Online optimization: Competing with dynamic comparators.
\newblock In {\em Artificial Intelligence and Statistics}, pages 398--406,
  2015.

\bibitem{kim2017online}
S.-J. Kim and G.~B. Giannakis.
\newblock An online convex optimization approach to real-time energy pricing
  for demand response.
\newblock {\em IEEE Transactions on Smart Grid}, 8(6):2784--2793, 2017.

\bibitem{lesage2018setpoint}
A.e Lesage-Landry and J.A. Taylor.
\newblock Setpoint tracking with partially observed loads.
\newblock {\em IEEE Transactions on Power Systems}, 33(5):5615 -- 5627, 2018.

\bibitem{mathieu2013state}
J.~L. Mathieu, S.~Koch, and D.~S. Callaway.
\newblock State estimation and control of electric loads to manage real-time
  energy imbalance.
\newblock {\em IEEE Transactions on Power Systems}, 28(1):430--440, 2013.

\bibitem{mohri2016accelerating}
M.~Mohri and S.~Yang.
\newblock Accelerating online convex optimization via adaptive prediction.
\newblock In {\em Artificial Intelligence and Statistics}, pages 848--856,
  2016.

\bibitem{mokhtari2016online}
A.~Mokhtari, S.~Shahrampour, A.~Jadbabaie, and A.~Ribeiro.
\newblock Online optimization in dynamic environments: Improved regret rates
  for strongly convex problems.
\newblock In {\em Decision and Control (CDC), 2016 IEEE 55th Conference on},
  pages 7195--7201. IEEE, 2016.

\bibitem{nesterov1998introductory}
Y.~Nesterov.
\newblock Introductory lectures on convex programming volume i: Basic course.
\newblock {\em Lecture notes}, 1998.

\bibitem{palensky2011demand}
P.~Palensky and D.~Dietrich.
\newblock Demand side management: Demand response, intelligent energy systems,
  and smart loads.
\newblock {\em IEEE transactions on industrial informatics}, 7(3):381--388,
  2011.

\bibitem{rakhlin2013online}
A.~Rakhlin and K.~Sridharan.
\newblock Online learning with predictable sequences.
\newblock In {\em Proceedings of the 26th Annual Conference on Learning Theory
  (COLT)}, pages 1--27, 2013.

\bibitem{rakhlin2013optimization}
S.~Rakhlin and K.~Sridharan.
\newblock Optimization, learning, and games with predictable sequences.
\newblock In {\em Advances in Neural Information Processing Systems}, pages
  3066--3074, 2013.

\bibitem{schmidt2011convergence}
M.~Schmidt, N.~L. Roux, and F.~R. Bach.
\newblock Convergence rates of inexact proximal-gradient methods for convex
  optimization.
\newblock In {\em Advances in neural information processing systems}, pages
  1458--1466, 2011.

\bibitem{shalev2012online}
S.~Shalev-Shwartz.
\newblock Online learning and online convex optimization.
\newblock {\em Foundations and Trends{\textregistered} in Machine Learning},
  4(2):107--194, 2012.

\bibitem{steinhardt2014adaptivity}
J.~Steinhardt and P.~Liang.
\newblock Adaptivity and optimism: An improved exponentiated gradient
  algorithm.
\newblock In {\em International Conference on Machine Learning}, pages
  1593--1601, 2014.

\bibitem{taylor2016power}
J.~A. Taylor, S.~V. Dhople, and D.~S. Callaway.
\newblock Power systems without fuel.
\newblock {\em Renewable and Sustainable Energy Reviews}, 57:1322--1336, 2016.

\bibitem{wright1999numerical}
S.~Wright and J.~Nocedal.
\newblock Numerical optimization.
\newblock {\em Springer Science}, 35(67-68), 1999.

\bibitem{yang2016optimistic}
S.~Yang and M.~Mohri.
\newblock Optimistic bandit convex optimization.
\newblock In {\em Advances in Neural Information Processing Systems}, pages
  2297--2305, 2016.

\bibitem{zinkevich2003online}
M.~Zinkevich.
\newblock Online convex programming and generalized infinitesimal gradient
  ascent.
\newblock In {\em Proceedings of the 20th International Conference on Machine
  Learning (ICML-03)}, pages 928--936, 2003.

\end{thebibliography}
\end{document}